\documentclass[12pt]{article} 
\usepackage[bookmarks=false, backref=page]{hyperref}
\usepackage[left=2.00cm, right=2.00cm, top=2.50cm, bottom=2.50cm]{geometry}
\usepackage{setspace}
\usepackage{abstract}
\usepackage{cite}
\usepackage{amsthm}
\usepackage{amsmath}
\usepackage{amsfonts}
\usepackage{amssymb}
\usepackage{mathrsfs}
\usepackage{color}
\usepackage{graphicx}
\usepackage{multirow}
\usepackage{booktabs}
\usepackage{indentfirst}
\usepackage{arydshln}
\usepackage{hyperref}
\allowdisplaybreaks[4]
\numberwithin{equation}{section}
\usepackage{authblk}
\makeatletter
\newcommand{\rmnum}[1]{\romannumeral #1}
\newcommand{\Rmnum}[1]{\expandafter\@slowromancap\romannumeral #1@}
\makeatother

\newtheorem{theorem}{Theorem}[section]

\newtheorem{proposition}[theorem]{Proposition}

    \newcommand{\BQ}{{\mathbb {Q}}} \newcommand{\BR}{{\mathbb {R}}}

\newcommand{\Res}{{\mathrm{Res}}}

    \DeclareFontFamily{U}{wncy}{}
\DeclareFontShape{U}{wncy}{m}{n}{<->wncyr10}{}
\DeclareSymbolFont{mcy}{U}{wncy}{m}{n}
\DeclareMathSymbol{\Sha}{\mathord}{mcy}{"58}

\author{Tong Wei and Shuai Zhai}
\title{Binary quadratic forms and elliptic curves with analytic rank one}

\begin{document}
\date{}
\maketitle

\vspace{-1cm}

\begin{abstract}
Given an elliptic curve with Weierstrass equation $y^2=f(x)$, and a positive definite binary quadratic form $Q(u, v)$. We show that there are infinitely many $d$ in the set represented by the quadratic forms in the genus of $Q$ such that the twisted elliptic curve $dy^2=f(x)$ has analytic rank one.
\end{abstract}

\section{Introduction}

Let $E: y^2 = f(x)$ be an elliptic curve defined over $\mathbb{Q}$ with conductor $q$. Let $L(E, s)$ be the complex $L$-series of $E$. For each discriminant $d$ of a quadratic extension of $\BQ$, we write $E^{(d)}$ for the twist of $E$ by this quadratic extension. Let $Q(u,v)$ be a positive definite binary quadratic form of discriminant $D$ over $\mathbb{Q}$,  and let $r_Q(d)$ denote the number of representations of the integer $d$ by the quadratic forms in the genus of $Q$. We define a family of elliptic curves given by
\begin{equation}\label{eq-Euv}
E_{u,v} : Q(u,v) y^2 = f(x),\quad u,v\in\mathbb{Z},\quad \gcd(u,v)=1.
\end{equation}

Previously, linear families of quadratic twists were studied by Goldfeld--Hoffstein--Patterson \cite{Gold81} for elliptic curves with complex multiplication, and by Murty--Murty \cite{MM} and Bump--Friedberg--Hoffstein \cite{Bu90} for general elliptic curves. For simplicity, throughout this paper, we assume that $(D,q)=1$. If $q$ is a square, we further assume that there exists a pair $(u_0,v_0)$ such that the elliptic curve $E_{u_0,v_0}$ has root number $-1$. Munshi \cite{Mu3,Mu4} made the first substantial progress on nonlinear families of quadratic twists of elliptic curves with complex multiplication. He proved that there exist infinitely many coprime integer pairs $(u,v)$ such that the elliptic curve $E_{u,v}$ has analytic rank one.

In this paper, we extend Munshi's work \cite{Mu4} from elliptic curves with complex multiplication to arbitrary elliptic curves. Under the above assumptions, we prove the following theorem.
\begin{theorem}\label{MainThm}
There exist infinitely many coprime integer pairs $(u,v)$ such that the elliptic curve $E_{u,v}$ has analytic rank one.
\end{theorem}

As presenting in \cite{Mu4}, there are several applications corresponding to the arithmetic of elliptic curves. For example, any positive definite binary quadratic form of odd discriminant properly represents infinitely many congruent numbers. While our result can be widely applied on the Diophantine problems associated with elliptic curves without complex multiplication, such as the generalised congruent number problem (see, for example, Top--Yui  \cite{TY}). Theorem \ref{MainThm} will follow by establishing the following weighted mean value theorem.
\begin{theorem}\label{thm-1.1}
We have
\begin{equation*}
{\sum_{(d, 2qD)=1}}^*r_Q(d)L'(1, E^{(d)}) F(d/X) = \alpha X\log X+O(X(\log X)^{\frac{1}{2}}(\log\log X)^3),
\end{equation*}
where the big-$O$ is depending on $D$ and $q$, $\sum^{*}$ denote a sum over square-free integers, and $F$ is a nonnegative smooth function with compact support. The leading coefficient $\alpha \neq 0$ when either of the following conditions hold: $(\rmnum{1})$ $q$ is not a square, $(\rmnum{2})$ the root number of $E$ is $-1$.
\end{theorem}

Our proof builds on the celebrated work of Munshi \cite{Mu4} together with Li's recent breakthrough on second moments \cite{Li}. This breakthrough has become a fundamental ingredient in a number of subsequent developments, including \cite{H,KMSS,Zh,J,W-Zh,H26}.

\medskip
{\bf Acknowledgements}. The authors would like to thank Yongxiao Lin and Jianya Liu for their helpful comments, and Zhining Wei for valuable discussions during his visit to Shandong University in June 2026, which helped simplify the proof of Proposition 3.2.

\section{Preliminaries}
Using genus theory, a square-free positive integer $d$ is represented by a form in the genus of $Q$ if and only if
\begin{equation}\label{eq.1.1}
r_Q(d)=\prod_{i=1}^{h}(\psi_i(d)+\varepsilon_i)\prod_{p\mid d}(1+\chi_D(p))
\end{equation}
is non-zero, where $h=\omega(D)$ or $\omega(D)\pm 1$, $\psi_i$ is the genus character, $\varepsilon_i=\pm 1$ and $\prod_i \varepsilon_i = 1$. (Here $\omega(D)$ denotes the number of distinct prime factors of $D$.) We set
\begin{equation}\label{eq.1.2}
R_D(d)=\prod_{p\mid d}(1+\chi_D(p)),
\end{equation}
and note that it is a multiplicative function on the set of square-free integers. For a more detailed discussion of genus theory, see \cite[Section 2, 3]{Mu4} or \cite[p.505--p.508]{I}.

Let $E/\mathbb{Q}$ be any elliptic curve of conductor $q$. There exists a primitive cusp form with level $q$ and weight $2$ as follow
\begin{equation*}
f(z)=\sum_{n=1}^{\infty}\lambda_{f}(n)n^{\frac{1}{2}}e(nz)\in S_{2}(\Gamma_0(q)),
\end{equation*}
such that
\begin{equation*}
L(s+1/2, E)=L(s, f)=\sum_{n=1}^{\infty}\lambda_{f}(n)n^{-s}
\end{equation*}
with $\lambda_f(1)=1$ and $f$ has been normalized so that the Deligne's bound gives $|\lambda_f(n)|\leq\tau(n)$ for all $n\geq1$, where $\tau(n)$ is the divisor function. We study the family of twists of $f$ by quadratic characters. Let $d$ denote a fundamental discriminant, and $\chi_d(\cdot)=\left(\frac{d}{\cdot}\right)$ denote the primitive quadratic character of conductor $\lvert d\rvert$. Then $f\otimes\chi_d$ is a primitive Hecke eigenform of level $q\lvert d\rvert^{2}$ and the twisted $L$-function is defined for $\Re(s)>1$ by
\begin{equation*}
L\left(s, f \otimes \chi_d\right)=\sum_{n=1}^{\infty}\frac{\lambda_f(n)\chi_d(n)}{n^s}=\prod_{p\nmid qd}\left(1-\frac{\lambda_f(p)\chi_{d}(p)}{p^s}+\frac{1}{p^{2s}}\right)^{-1}
\prod_{p\mid q}\left(1-\frac{\lambda_f(p)\chi_{d}(p)}{p^s}\right)^{-1}.
\end{equation*}
This satisfies the functional equation
\begin{equation}
\Lambda\left(s, f \otimes \chi_d\right)=\varepsilon_f\chi_d(-q^*)\Lambda\left(1-s, f \otimes \chi_d\right),
\end{equation}
where $q^*$ is the square-free part of $q$, $\varepsilon_f$ is the sign of the functional equation of $L(s, f)$ (or the root number of the elliptic curve $E$), and the completed $L$-function is given by
\begin{equation}
\Lambda\left(s, f \otimes \chi_d\right)=\left(\frac{\sqrt{q}\lvert d\rvert}{2 \pi}\right)^{s} \Gamma\left(s+1/2\right) L\left(s, f \otimes \chi_d\right).
\end{equation}
To evaluate the $L$-function inside the critical strip we can use the functional equation to get a rapidly decaying series expansion called the approximate functional
equation. For the first derivative at the center, if $\varepsilon_f\chi_d(-q^*)=-1$ we have
\begin{equation}\label{eq.L}
L'\left(1/2, f \otimes \chi_d\right)=(1-\varepsilon_f\chi_d(-q^*))\sum_{n\geq 1}\frac{\lambda_f(n)\chi_d(n)}{\sqrt{n}}W\left(\frac{n}{|d|}\right),
\end{equation}
where the smooth function $W$ is given by
\begin{equation*}
W(y)=\frac{1}{2\pi i}\int_{(3)}\frac{\Gamma(s+1)}{(2\pi/\sqrt{q})^s}y^{-s}\frac{ds}{s^2}
\end{equation*}
This function decays rapidly as $y$ gets larger, we have
$$ W(y)\ll_M y^{-M} $$
for all $M\geq 1$. This implies that in \eqref{eq.L} the terms with $n \gg d$ make a negligibly small contribution, in other words the first $\asymp d$ terms of the sum give a very good approximation for the central value. When $y$ is small, we have 
\[
W(y) = (-\log y) + v + O(y)
\]
for some constant $v$. This fact will play an important role in our analysis.

\section{Proof of Theorem 1.2}
For simplicity, we will assume that $(D, q)=1$. In fact, our method is valid for the general cases, and can be worked out without this coprimality assumption. Let $Q$ be the principal quadratic form of discriminant $D$. We are interested in the sum
\begin{equation}\label{eq.1.3}
{\sum_{(d, 2qD)=1}}^*(1+\chi_{-4}(d))r(d)L'(1/2, f\otimes \chi_D \chi_d) F(d/X).
\end{equation}
The weight $r(d)$ is given by
\begin{equation}\label{eq.1.4}
r(d)=\prod_{p\mid D}\left(1+\left(\frac{p}{d}\right)\right)\times R_D(d).
\end{equation}
which is related (via reciprocity) with $r_Q(d)$ as defined in \eqref{eq.1.1} for the principal form $Q$. In fact, in any case $(1+\chi_{-4}(d))r(d)\neq0$ implies that $r_Q(d)\neq0$. We are imposing a condition so that only square-free odd integers $d\equiv 1\pmod 4$ contribute to the sum \eqref{eq.1.3}. This is necessary for using reciprocity. We are also throwing in an auxiliary character $\chi_D$ to simplify the computation of the main term.

With hindsight, we set $Y=\frac{X}{(\log X)^{100}}$. Further let
\begin{equation}\label{eq-A,B}
\begin{aligned}
\mathcal{A}(d)&=(1-\varepsilon_f\chi_d(-q^*))\sum_{n\geq 1}\frac{\lambda_f(n)\chi_d(n)}{\sqrt{n}}\chi_D(n)W\left(\frac{n}{Y}\right),\\
\mathcal{B}(d)&=L'(1/2, f\otimes \chi_D \chi_d)-\mathcal{A}(d).
\end{aligned}
\end{equation}
To prove Theorem \ref{thm-1.1}, we only need to verify the following two propositions; their proofs appear in Sections \ref{pro-1.1} and \ref{pro-1.2}, respectively.
\begin{proposition}\label{pro-1.1}
With notation as above, we have
\begin{equation*}
\begin{aligned}
&\quad{\sum_{(d, 2qD)=1}}^*(1+\chi_{-4}(d))r(d)\mathcal{A}(d)F(d/X)\\
&=2C\left(1-\varepsilon_f\prod_{p\mid q^*}\frac{\lambda_f(p)\chi_D(p)}{\sqrt{p}}\left(1+\frac{1}{p}\right)^{-1}\right)X\log X+O(X\log\log X),
\end{aligned}
\end{equation*}
for some constant $C$, the big-$O$ is depending on $D$ and $q$.
\end{proposition}
\begin{proposition}\label{pro-1.2}
With notation as above, we have
\begin{equation*}
{\sum_{(d, 2qD)=1}}^*(1+\chi_{-4}(d))r(d)\mathcal{B}(d)F(d/X)\ll X(\log X)^{\frac{1}{2}}(\log\log X)^3,
\end{equation*}
the big-$O$ is depending on $D$ and $q$.
\end{proposition}
\subsection{Proof of Proposition 3.1}
Using the approximate functional equation \eqref{eq.L} and \eqref{eq.1.4}, observe that whenever $d \equiv 1 \pmod{4}$, then $\chi_d(n) = \chi_n(d)$, we thus write
\begin{equation}\label{eq.S}
\begin{aligned}
S&={\sum_{(d, 2qD)=1}}^*(1+\chi_{-4}(d))r(d)\mathcal{A}(d)F(d/X)\\
&={\sum_{(d, 2qD)=1}}^*\left(1+\left(\frac{-1}{d}\right)\right)\left(1-\varepsilon_f\left(\frac{q^*}{d}\right)\right)\prod_{p\mid D}\left(1+\left(\frac{p}{d}\right)\right) R_D(d)\\
&\quad\times\sum_{n=1}^{+\infty}\frac{\lambda_f(n)\chi_D(n)}{\sqrt{n}}\left(\frac{n}{d}\right)W\left(\frac{n}{Y}\right)F\left(\frac{d}{X}\right).
\end{aligned}
\end{equation}
Hence, the computation of $S$ boils down to computing sums of the form
\begin{equation}\label{eq.T}
T={\sum_{(d, 2qD)=1}}^*R_D(d)\left(\frac{a}{d}\right)\sum_{n=1}^{+\infty}\frac{\lambda_f(n)\chi_D(n)}{\sqrt{n}}\left(\frac{n}{d}\right)W\left(\frac{n}{Y}\right)F\left(\frac{d}{X}\right),
\end{equation}
for a given square-free integer $a$ (which is $\delta\theta$ where $\delta\mid D$ and $\theta=\pm1, \pm q^*$). Interchanging the order of summation and applying the Mellin inverse transform of $F$, we obtain the following expression for \eqref{eq.T},
\begin{equation}\label{eq.T1}
\begin{aligned}
T&=\frac{1}{2\pi i}\int_{(c)}\widetilde{F}(s)X^s\sum_{n\geq 1}\frac{\lambda_f(n)\chi_D(n)}{\sqrt{n}}W\left(\frac{n}{Y}\right)
{\sum_{(d, 2qD)=1}}^*\frac{R_D(d)}{d^s}\left(\frac{an}{d}\right)ds\\
&=\frac{1}{2\pi i}\int_{(c)}\widetilde{F}(s)X^s\sum_{n\geq 1}\frac{\lambda_f(n)\chi_D(n)}{\sqrt{n}}W\left(\frac{n}{Y}\right)
\frac{L(s, \chi_{4an})L(s, \chi_{4Dan})}{\zeta(2s)^{2}L(2s, \chi_D)}D_{n}(s)ds.
\end{aligned}
\end{equation}
The Dirichlet series $D_{n}(s)$ converges absolutely for $\Re(s)=\sigma>1/3$, and in the domain $\sigma\geq 1/2$ it satisfies the bound $|D_{n}(s)|\ll \frac{n}{\varphi(n)}$, where the implied constant depends only on $q, D$. We shift the contour of integration to the critical line $\sigma=1/2$ and in the process we pick up residue at the only possible pole at $s=1$. This gives $T=M+E$, where $M$ comes from the residue at $s=1$ and it contributes to the main term. The term $E$ is given by the integral over the central line $\sigma=1/2$ and it contributes to the error term. For $E$, we observe that the terms with $n>Y^{1+\varepsilon}$ for any $\varepsilon>0$ contribute a negligible quantity compared to the other terms, and hence can be ignored. So it follows that
\begin{equation*}
E\ll\int_{(1/2)}|\widetilde{F}(s)D_{n}(s)|\frac{\sqrt{X}}{|\zeta(2s)^{2}L(2s, \chi_D)|}\sum_{n\ll Y^{1+\varepsilon}}\frac{|\lambda_f(n)|}{\sqrt{n}}\left|W\left(\frac{n}{Y}\right)\right|
|L(s, \chi_{4an})L(s, \chi_{4Dan})||ds|.
\end{equation*}
\subsubsection{The main term}
We analyse the main term
\begin{equation*}
M=\sum_{n\geq 1}\frac{\lambda_f(n)\chi_D(n)}{\sqrt{n}}W\left(\frac{n}{Y}\right)\Res_{s=1}\left\{\widetilde{F}(s)X^s
\frac{L(s, \chi_{4an})L(s, \chi_{4Dan})}{\zeta(2s)^{2}L(2s, \chi_D)}D_{n}(s)\right\}.
\end{equation*}
There is a simple pole when either $an=\square$ or $Dan=\square$. Recall that  $a=\delta\theta$ where $\delta\mid D$ and $\theta=\pm1, \pm q^*$. Due to the presence of the character $\chi_D$ we have $(D, n)=1$. So $an=\square$ can only hold if $a=1, q^*$ and similarly $Dan=\square$ can only hold if $a=D, Dq^*$. Thus, we only focus on these four cases in the evaluation of the main term.

To compute the residue we go back to the expression \eqref{eq.T1} and observe that
\begin{equation*}
\frac{L(s, \chi_{4an})L(s, \chi_{4Dan})}{\zeta(2s)^{2}L(2s, \chi_D)}D_{n}(s)=\prod_{p\nmid 2qD}
\left(1+\frac{1}{p^s}\left(\frac{an}{p}\right)+\frac{1}{p^s}\left(\frac{Dan}{p}\right)\right).
\end{equation*}
Then when $an=\square$ or $Dan=\square$, we have
\begin{equation}\label{eq-res}
\Res_{s=1}\left\{\frac{L(s, \chi_{4an})L(s, \chi_{4Dan})}{\zeta(2s)^{2}L(2s, \chi_D)}D_{n}(s)\right\}=C_{q, D}\prod_{\substack{p\mid n \\ p\nmid 2qD}}
\left(1+\frac{1}{p}+\frac{1}{p}\left(\frac{D}{p}\right)\right)^{-1},
\end{equation}
for some constant $C_{q, D}\neq 0$. So we obtain that
\begin{equation*}
M=C_{q, D}\widetilde{F}(1) X {\sum_{n\geq 1}}^{\flat}\frac{\lambda_f(n)\chi_D(n)}{\sqrt{n}}W\left(\frac{n}{Y}\right)\prod_{\substack{p\mid n \\ p\nmid 2qD}}
\left(1+\frac{1}{p}+\frac{1}{p}\left(\frac{D}{p}\right)\right)^{-1},
\end{equation*}
where $\flat$ indicates that either $an=\square$ or $Dan=\square$. It follows from the definition of \(W\) that
\begin{equation*}
=C_{q, D}\widetilde{F}(1) X \frac{1}{2\pi i}\int_{(3)}Y^s\frac{\Gamma(s+1)}{(2\pi/\sqrt{q})^s}D(s)\frac{ds}{s^2}
\end{equation*}
the Dirichlet series $D(s)$ is as follow
\begin{equation}\label{eq-D_1}
D(s)={\sum_{n\geq 1}}^{\flat}\frac{\lambda_f(n)\chi_D(n)}{n^{\frac{1}{2}+s}}\prod_{\substack{p\mid n \\ p\nmid 2qD}}
\left(1+\frac{1}{p}+\frac{1}{p}\left(\frac{D}{p}\right)\right)^{-1}.
\end{equation}
For convenience, we write $(1+\tfrac{1}{p}+\tfrac{1}{p}(\tfrac{D}{p}))^{-1}=l_p$. Using the multiplicativity of the Fourier coefficient we get that when $a = 1, D$,
\begin{equation*}
\begin{aligned}
D(s)&=\prod_{p\nmid 2qD}\left\{1+\frac{\lambda_f(p^2) l_p}{p^{1+2s}}+\sum_{v\geq 2}\frac{\lambda_f(p^{2v})l_p}{p^{2v(\frac{1}{2}+s)}}\right\}
\prod_{p\mid 2q}\sum_{v\geq 0}\frac{\lambda_f(p^{2v})}{p^{2v(\frac{1}{2}+s)}},
\end{aligned}
\end{equation*}
when $a = q^*, Dq^*$,
\begin{equation*}
D(s)=\prod_{p\nmid 2qD}\left\{1+\frac{\lambda_f(p^2) l_p}{p^{1+2s}}+\sum_{v\geq 2}\frac{\lambda_f(p^{2v})l_p}{p^{2v(\frac{1}{2}+s)}}\right\}
\prod_{p\mid 2q}\sum_{v\geq 0}\frac{\lambda_f(p^{2v})}{p^{2v(\frac{1}{2}+s)}}
\times \prod_{p\mid q^*}\frac{\lambda_f(p)\chi_D(p)}{p^{\frac{1}{2}+s}}\left(1+\frac{1}{p^{1+2s}}\right)^{-1}.
\end{equation*}

From the above, it can be seen that
\begin{equation}\label{eq-D_2}
D(s)=L(1+2s, \operatorname{sym}^2 f)D^*(s),
\end{equation}
here $D^*(s)$ is absolutely convergent and has a Euler product expression in $\Re(s)\geq -1/4+\varepsilon$. We now shift the contour of integration in $s$ to $\Re(s)=-1/4+\varepsilon$. We use the standard convexity estimate$L^{(j)}(\sigma+it, \operatorname{sym}^2 f)\ll(1+|t|)^{\frac{3}{2}(1-\sigma)}(\log(2+|t|))^{j+1}$ for $0\leq\sigma\leq1$, it is evident that the contribution of the contour is $\ll_{q, D} X^{\frac{3}{4}+\varepsilon}$.
The contribution of the residue at $s=0$ is
\begin{equation*}
C_{q, D}\widetilde{F}(1)D(0) X\log Y +C'X
\end{equation*}
for some constant $C'$. Hence, the leading term of $S$ as defined in \eqref{eq.S} is given by
\begin{equation*}
2C\left(1-\varepsilon_f\prod_{p\mid q^*}\frac{\lambda_f(p)\chi_D(p)}{\sqrt{p}}\left(1+\frac{1}{p}\right)^{-1}\right)X\log X,
\end{equation*}
here $C$ is a constant.
\subsubsection{The error term}\label{sub-1.2}
Let us recall that
\begin{equation*}
E\ll\int_{(1/2)}|\widetilde{F}(s)D_{n}(s)|\frac{\sqrt{X}}{|\zeta(2s)^{2}L(2s, \chi_D)|}\sum_{n\ll Y^{1+\varepsilon}}\frac{|\lambda_f(n)|}{\sqrt{n}}\left|W\left(\frac{n}{Y}\right)\right|
|L(s, \chi_{4an})L(s, \chi_{4Dan})||ds|.
\end{equation*}
We will break up the inner sum into dyadic blocks and in the block $N<n\leq 2N$ we will use the bound $W(n/Y)\ll H(N)$, where
\begin{equation}\label{eq-H(N)}
H(N)=
\begin{cases}
\log(Y/N)+1, & \text{if } N\leq Y,\\
Y/N, & \text{if } Y<N\ll Y^{1+\varepsilon}.
\end{cases}
\end{equation}
Then using inequality $2xy\leq x^2+y^2$, and using the trivial bound for zeta function and Dirichlet $L$-function on the vertical line, and $\max_{n\ll Y^{1+\varepsilon}}\sum_{p\mid n}\frac{1}{p}\asymp\log\log\log Y$, $\widetilde{F}(s)\ll(1+|s|)^{-A}$ for any $A>0$, we get that
\begin{equation*}
E\ll\sqrt{X}\log\log Y\sum_{\text{dyadic blocks}}H(N)N^{-\frac{1}{2}}\int_{\BR}U(N, t)(1+|t|)^{-7}dt,
\end{equation*}
where the sum is over $\log Y$ many dyadic blocks, and
\begin{equation*}
U(N, t)=\sum_{N<n\leq 2N}|\lambda_f(n)||L(s, \chi_{4bn})|^2,
\end{equation*}
note that $b=a$ or $b=Da$. Following the method outlined in Sections 6-8 of Munshi's paper \cite{Mu4}, we obtain from \cite[Lemma 5]{Mu4} that
\begin{equation*}
U(N, t)\ll\delta(t)N(\log N)^{\frac{3}{2}-\theta},
\end{equation*}
where $\delta(t)$ grow like $|t|^2$ as $|t|\rightarrow\infty$ and near $t=0$ it blows up like $|t|^{-2}$, $\theta=1/2-2^{\frac{1}{4}}/3$. This gives us
\begin{equation*}
\begin{aligned}
E&\ll\sqrt{X}\log\log Y\sum_{\text{dyadic blocks}}H(N)N^{\frac{1}{2}}(\log N)^{\frac{3}{2}-\theta}\int_{\BR}\delta(t)(1+|t|)^{-7}dt\\
&\ll \sqrt{X}\log\log Y(\log Y)^{\frac{3}{2}-\theta}\sum_{\text{dyadic blocks}}H(N)N^{\frac{1}{2}}
\end{aligned}
\end{equation*}
We add up the contribution of all the dyadic blocks using the bound in \eqref{eq-H(N)}, and taking $Y=X/(\log X)^{100}$, we have
\begin{equation*}
E\ll X^{\frac{1}{2}}Y^{\frac{1}{2}}(\log\log Y)^3(\log Y)^{\frac{3}{2}-\theta}\ll X.
\end{equation*}
This completes the proof of Proposition \ref{pro-1.1}.

\subsection{Proof of Proposition 3.2}

From the approximate functional equation \eqref{eq.L} and \eqref{eq.1.4}, \eqref{eq-A,B}, we see that
\begin{equation*}
\begin{aligned}
&\quad{\sum_{(d, 2qD)=1}}^*(1+\chi_{-4}(d))r(d)\mathcal{B}(d)F(d/X)\\
&={\sum_{(d, 2qD)=1}}^*(1+\chi_{-4}(d))(1-\varepsilon_f\chi_d(-q^*))r(d)\\
&\quad\times\sum_{n\geq 1}\frac{\lambda_f(n)\chi_d(n)}{\sqrt{n}}\chi_D(n)\left\{W\left(\frac{n}{|d|}\right)-W\left(\frac{n}{Y}\right)\right\}F\left(\frac{d}{X}\right).
\end{aligned}
\end{equation*}
There exists a smooth real-valued function $G$ compactly supported on $[3/4, 2]$ which satisfies
\begin{equation}\label{eq.G}
\begin{aligned}
G(x)=1 &\text{ for all $x\in[1, 3/2]$},\\
G(x)+G(x/2)=1 &\text{ for all $x\in[1,3]$}.
\end{aligned}
\end{equation}
For more detailed information, readers can refer to Warner's book\cite[Theorem 1.11 and Corollary]{Warner} or \cite[p709]{Li}.
Let $H$ be a positive integer. The function $F(x)$ defined by
\begin{equation}\label{eq.F}
F(x)=G(x)+G(x/2)+\ldots+G\left(x/2^H\right)
\end{equation}
satisfies $F(x)=1$ for all $x\in [1, 3\cdot 2^{H-1}]$ and support on $[3/4, 2^{H+1}]$. We have the smooth partition of unity
\begin{equation}
\sum_{H=0}^{\infty} G\left(\frac{x}{2^H}\right)=1
\end{equation}
for all $x\in [1,\infty)$.

A dyadic partition of unity using the functions $G$, we obtain
\begin{equation}\label{eq-pro-1.2}
\begin{aligned}
&\quad{\sum_{(d, 2qD)=1}}^*(1+\chi_{-4}(d))r(d)\mathcal{B}(d)F(d/X)\\
&={\sum_{N}}^{d}{\sum_{(d, 2qD)=1}}^*(1+\chi_{-4}(d))(1-\varepsilon_f\chi_d(-q^*))r(d)\\
&\quad\times\sum_{n\geq 1}
\frac{\lambda_f(n)\chi_d(n)}{\sqrt{n}}\chi_D(n)\left\{W\left(\frac{n}{|d|}\right)-W\left(\frac{n}{Y}\right)\right\}F\left(\frac{d}{X}\right)G\left(\frac{n}{N}\right),
\end{aligned}
\end{equation}
where $\sum^{d}$ denotes the sum over $N=2^{H}$ for integer $H\geq0$. We note that
\begin{equation}\label{eq-r}
{\sum_{(d, 2qD)=1}}^*(1+\chi_{-4}(d))^2(1-\varepsilon_f\chi_d(-q^*))^2 r(d)^2 F\left(\frac{d}{X}\right)\ll_D \sum_{d=1}^{\infty}R_{D}(d)^2F\left(\frac{d}{X}\right).
\end{equation}
By the inverse Mellin transform for $F$, we obtain
\begin{equation*}
\sum_{d=1}^{\infty}R_{D}(d)^2F\left(\frac{d}{X}\right)=\frac{1}{2\pi i}\int_{(3)}\widetilde{F}(s)X^s L_D(s)ds,
\end{equation*}
where
$$\sum_{d=1}^{\infty}\frac{R_{D}(d)^2}{d^s}=\zeta(s)^2L(s, \chi_D)^2Z(s),$$
\(Z\) converges for \(\Re(s) > 1/2\). We shift the contour to $\Re(s)=3/4$. Noticing that we encounter a second order pole at $s=1$ for $L_D(s)$, then we have
\begin{equation}\label{eq-r1}
\sum_{d=1}^{\infty}R_{D}(d)^2F\left(\frac{d}{X}\right)\ll_D X\log X.
\end{equation}

We split \eqref{eq-pro-1.2} into three cases: $N>X$, $Y<N\leq X$ and $N\leq Y$. For $N>X$, applying the Cauchy--Schwarz inequality for the summation in $d$, combined with \cite[Section 7.5]{Zh}, \eqref{eq-r}, \eqref{eq-r1}, we see that the contribution from this part is $\ll X(\log X)^{\frac{1}{2}}$. Similarly, by the Cauchy--Schwarz inequality and \cite[Lemma 7.1]{Zh}, \eqref{eq-r}, \eqref{eq-r1}, we have that the contribution for $Y<N\leq X$ satisfies $\ll X(\log X)^{\frac{1}{2}}(\log\log X)^3$.

When $N\leq Y$, using inverse Mellin transform for $F$, our goal is to establish an upper bound of the following form
\begin{equation*}
\begin{aligned}
&\quad{\sum_{N\leq Y}}^{d}\frac{1}{2\pi i}\int_{(c)}\widetilde{F}(s)X^{s}\sum_{n\geq 1}\frac{\lambda_f(n)\chi_D(n)}{\sqrt{n}}G\left(\frac{n}{N}\right)\\
&\times{\sum_{(d, 2qD)=1}}^*(1+\chi_{-4}(d))(1-\varepsilon_f\chi_d(-q^*))r(d)\chi_d(n)d^{-s}\left\{W\left(\frac{n}{|d|}\right)-W\left(\frac{n}{Y}\right)\right\}ds.
\end{aligned}
\end{equation*}
The summation over \(d\) follows from \eqref{eq.S} and \eqref{eq.T}. Combined with the definition of $W$, we have
\begin{equation*}
\begin{aligned}
&{\sum_{(d, 2qD)=1}}^*(1+\chi_{-4}(d))(1-\varepsilon_f\chi_d(-q^*))r(d)\chi_d(n)d^{-s}\left\{W\left(\frac{n}{|d|}\right)-W\left(\frac{n}{Y}\right)\right\}\\
&=\sum_{\text{finite sum}}\frac{c' \varepsilon_f}{2\pi i}\int_{(c_1)}\frac{\Gamma(w+1)}{(2\pi/\sqrt{q})^w}n^{-w}\\
&\quad\times\left\{\frac{L(s-w, \chi_{4an})L(s-w, \chi_{4Dan})}{\zeta(2s-2w)^{2}L(2s-2w, \chi_D)}D_{n}(s-w)-\frac{L(s, \chi_{4an})L(s, \chi_{4Dan})}{\zeta(2s)^{2}L(2s, \chi_D)}D_{n}(s)Y^w\right\}\frac{dw}{w^2},
\end{aligned}
\end{equation*}
where $c'$ is a constant. We apply the change of variable $s\mapsto s+w$ on the first part of the above expression. It then follows that the contribution from the range $N\leq Y$ is
\begin{equation*}
\begin{aligned}
&\quad{\sum_{N\leq Y}}^{d}\sum_{\text{finite sum}}\frac{c' \varepsilon_f}{(2\pi i)^2}\int_{(c_2)}\int_{(c_1)}\frac{\Gamma(w+1)}{(2\pi/\sqrt{q})^w}\sum_{n\geq 1}
\frac{\lambda_f(n)\chi_D(n)}{n^{\frac{1}{2}+w}}G\left(\frac{n}{N}\right)\\
&\times\left\{\widetilde{F}(s+w)X^{s+w}\frac{L(s, \chi_{4an})L(s, \chi_{4Dan})}{\zeta(2s)^{2}L(2s, \chi_D)}D_{n}(s)-\widetilde{F}(s)X^{s}Y^w \frac{L(s, \chi_{4an})L(s, \chi_{4Dan})}{\zeta(2s)^{2}L(2s, \chi_D)}D_{n}(s)\right\}\frac{dw}{w^2}ds
\end{aligned}
\end{equation*}
for $c_1, c_2>0$. We shift the contour of integration in $s$ to $\Re(s)=1/2$, which splits the above expression into two parts, denoted by \(R\) and \(I\), where
\begin{equation*}
\begin{aligned}
R&={\sum_{N\leq Y}}^{d}\sum_{\text{finite sum}}c' \varepsilon_f\frac{X}{2\pi i}\int_{(c_1)}\frac{\Gamma(w+1)}{(2\pi/\sqrt{q})^w}{\sum_{n\geq 1}}^{\flat}
\frac{\lambda_f(n)\chi_D(n)}{n^{\frac{1}{2}+w}}G\left(\frac{n}{N}\right)\left\{\widetilde{F}(1+w)X^{w}-\widetilde{F}(1)Y^w\right\}\\
&\quad\times\Res_{s=1}\left\{\frac{L(s, \chi_{4an})L(s, \chi_{4Dan})}{\zeta(2s)^{2}L(2s, \chi_D)}D_{n}(s)\right\}\frac{dw}{w^2}.
\end{aligned}
\end{equation*}
\begin{equation*}
\begin{aligned}
I&={\sum_{N\leq Y}}^{d}\sum_{\text{finite sum}}\frac{c' \varepsilon_f}{(2\pi i)^2}\int_{(1/2)}\int_{(c_1)}\frac{\Gamma(w+1)}{(2\pi/\sqrt{q})^w}\sum_{n\geq 1}
\frac{\lambda_f(n)\chi_D(n)}{n^{\frac{1}{2}+w}}G\left(\frac{n}{N}\right)\\
&\times\left\{\widetilde{F}(s+w)X^{s+w}\frac{L(s, \chi_{4an})L(s, \chi_{4Dan})}{\zeta(2s)^{2}L(2s, \chi_D)}D_{n}(s)-\widetilde{F}(s)X^{s}Y^w \frac{L(s, \chi_{4an})L(s, \chi_{4Dan})}{\zeta(2s)^{2}L(2s, \chi_D)}D_{n}(s)\right\}\frac{dw}{w^2}ds.
\end{aligned}
\end{equation*}
By the inverse Mellin transform for $G(x)$ and formula \eqref{eq-res}, \eqref{eq-D_1}, \eqref{eq-D_2}, we have that
\begin{equation*}
\begin{aligned}
R&={\sum_{N\leq Y}}^{d}\sum_{\text{finite sum}}c_{q, D} \varepsilon_f\frac{X}{2\pi i}\int_{(c_1)}\int_{(3)}\frac{\Gamma(w+1)}{(2\pi/\sqrt{q})^w}N^u \widetilde{G}(u)\\ &\quad\times\left\{\widetilde{F}(1+w)X^{w}-\widetilde{F}(1)Y^w\right\}L(1+2w+2u, \operatorname{sym}^2 f)D^*(w+u)du\frac{dw}{w^2}
\end{aligned}
\end{equation*}
for some constant $c_{q, D}$. We shift both the integration contours of \(u\) and \(w\) to \(\Re(u)=\Re(w)=-1/16\), and the contribution from the new integral region is \(\ll {\sum_{N\leq Y}}^{d}X(X^{-\frac{1}{16}}-Y^{-\frac{1}{16}})N^{-\frac{1}{16}}\ll X\). The contribution from the pole at \(w = 0\) is
\begin{equation*}
\begin{aligned}
&\quad{\sum_{N\leq Y}}^{d}\sum_{\text{finite sum}}c_{q, D} \varepsilon_f\frac{X}{2\pi i}\int_{(-1/16)}N^u \widetilde{G}(u)\\
&\times\Res_{w=0}\left\{\frac{\Gamma(w+1)}{(2\pi/\sqrt{q})^w}\left(\widetilde{F}(1+w)X^{w}-\widetilde{F}(1)Y^w\right)L(1+2w+2u, \operatorname{sym}^2 f)D^*(w+u)\frac{1}{w^2}\right\}du\\
&\ll {\sum_{N\leq Y}}^{d}XN^{-\frac{1}{16}}(\log X-\log Y)\ll X\log\log X.
\end{aligned}
\end{equation*}
We finally handle the term \(I\). Let $\eta=\frac{\log\log X}{\log X}$ and $C_{\eta}$ be the contour:
\begin{equation*}
C_{\eta}=\{it: |t|\geq \eta\}\cup\{\eta e^{i\vartheta}: \vartheta\in[-\pi/2, \pi/2]\}.
\end{equation*}
We decompose $I$ into two components: the linear segment and the semicircular arc. Similar to the function $W(x)$, for $x\leq 1$, $G(x)\ll 1$; $G(x)\ll x^{-A}$ for any $A>0$, if $x>1$.
Therefore, the contribution from $n>N^{1+\varepsilon}$ is negligible compared to that from small $n$. Consequently, the contribution from the linear segment is
\begin{equation*}
\begin{aligned}
&\ll{\sum_{N\leq Y}}^{d}X^{\frac{1}{2}}
\int_{-\infty}^{+\infty}|\widetilde{F}(1/2+it_1)|\sum_{n\ll N^{1+\varepsilon}}\frac{|\lambda_f(n)|}{n^{\frac{1}{2}}}\left|G\left(\frac{n}{N}\right)\right|\\
&\quad\times\left|\frac{L(1/2+it_1, \chi_{4an})L(1/2+it_1, \chi_{4Dan})}{\zeta(1+2it_1)^{2}L(1+2it_1, \chi_D)}D_{n}(1/2+it_1)\right|dt_1,
\end{aligned}
\end{equation*}
the contribution from the semicircular arc is
\begin{equation*}
\begin{aligned}
&\ll{\sum_{N\leq Y}}^{d}X^{\frac{1}{2}+\eta}\frac{1}{\eta^2}\int_{-\infty}^{+\infty}|\widetilde{F}(1/2+it)|\sum_{n\ll N^{1+\varepsilon}}
\frac{|\lambda_f(n)|}{n^{\frac{1}{2}}}\left|G\left(\frac{n}{N}\right)\right|\\
&\quad\times\left|\frac{L(1/2+it, \chi_{4an})L(1/2+it, \chi_{4Dan})}{\zeta(1+2it)^{2}L(1+2it, \chi_D)}D_{n}(1/2+it)\right|dt.
\end{aligned}
\end{equation*}
As shown in Section \ref{sub-1.2} ($G(x)\ll 1$, if $x\leq 1$; $G(x)\ll x^{-A}$ for any $A>0$, if $x>1$), $Y=X/(\log X)^{100}$, we have
\begin{equation*}
I\ll{\sum_{N\leq Y}}^{d}X^{\frac{1}{2}+\eta}\frac{1}{\eta^2}N^{\frac{1}{2}}(\log\log N)(\log N)^{\frac{3}{2}-\theta}\ll X.
\end{equation*}
This completes the proof of Proposition \ref{pro-1.2}.

\bigskip

\noindent Tong Wei, {\it Research Center for Mathematics and Interdisciplinary Sciences, Shandong University, Qingdao, Shandong, China.}

\smallskip
\noindent {\it E-mail:} 202421344@mail.sdu.edu.cn

\bigskip

\noindent Shuai Zhai, {\it Mathematical Research Center, Shandong University, Jinan, Shandong, China.}

\smallskip
\noindent {\it E-mail:} zhai@sdu.edu.cn

\end{document}